\documentclass[11pt]{amsart}
\usepackage{geometry}                % See geometry.pdf to learn the layout options. There are lots.
\usepackage{graphicx}
\usepackage{amssymb}
\usepackage{epstopdf}
\usepackage{ccaption} %% package for bilingual captions
\usepackage[english,french]{babel} 

\usepackage{subfig} % extra packages for figures and Taylor diagrams
\usepackage{color}
\usepackage{caption}
\usepackage{cancel}
\newenvironment{poliabstract}[1]
  {\renewcommand{\abstractname}{#1}\begin{abstract}}
  {\end{abstract}}

\newcommand {\R} {\mathbb {R}}

\newcommand {\Pj} {\mathbb {P}}

\theoremstyle{definition}

\newtheorem{defn}{Definition}[section]

\newtheorem{thm}{Theorem}[section]

\newtheorem{rem}[thm]{Remark}

\title{The title of my page}

\title[Une approche Semple \`a Goursat: le cas multi-drapeau]{A Semple Approach to Goursat: the multi-flag case.\\
Une approche Semple \`a Goursat: le cas multi-drapeau.}
\author{Alex L. Castro}
\email{alex.castro@mat.puc-rio.br}
\address{Departamento de Matem\'atica, Pont\'ificia Universidade Cat\'olica do Rio de Janeiro. Phone: +55-21-3527-1722}
\author{Wyatt C. Howard}
\email{whoward@ucsc.edu}
\thanks{We thank Gary Kennedy (Ohio State) and Richard Montgomery (UCSC) for many useful remarks, and Corey Shambron (UCSC) for attentively proofreading an earlier version of this paper. }
\address{Mathematics Department, University of California at Santa Cruz}
\date{\today}
\begin{document}

% ----------------------------------------------------------------

\maketitle
%% english abstract
\selectlanguage{english}
\begin{poliabstract}{Abstract} 

We consider the problem of classifying the orbits within a tower of fibrations with fibers diffeomorphic to projective planes and generalize the tower of fiber bundles due to J. Semple.  This tower was rediscovered by Montgomery and Zhitomirskii in the context of subriemannian geometry.
This tower admits a natural action of the diffeomorphism group of affine 3-space and these orbits correspond to classes of Goursat multi-flags. We demonstrate that it is possible to classify many of these orbits by elementary means by appealing to some basic tools in projective geometry, and the combinatorics of spatial curves.

% An abstract is a short summary of the main points in the project.  It
% does not mention background material or detail.

\end{poliabstract}

%%french abstract 

\selectlanguage{french}
\begin{poliabstract}{R\'esum\'e} 
Nous nous int\'eressons  au probl\`eme de classification des orbites dans
une tour de fibration, o\`u les fibres sont diff\'eomorphes \`a des plans
projectifs, g\'en\'eralisant ainsi les tours de fibr\'es projectives de J.
Semple. Cette tour a \'et\'e red\'ecouverte par R. Montgomery et M.
Zhitomirskii dans le contexte de la g\'eom\'etrie sub-riemannienne. Cette
tour admet une action naturelle du groupe des diff\'eomorphismes de
l'espace affine de dimension 3 et ces orbites correspondent \`a des
classes de multi-drapeaux de Goursat. Nous d\'emontrons quÕil est possible
de classifier un grand nombre de ces orbites de mani\'ere \'el\'ementaire en
utilisant des outils classiques de g\'eom\'etrie projective et la
combinatoire des courbes spatiales.
\end{poliabstract}

\vspace{.2in}


\noindent Secondary: $58A17$, $58K50$.
\vspace{.1in}

\vspace{.2in}

%%%%%%%%%%%%%%%%%%%%%%%%%%%%%%%%%%%%%%%%%%%%%%%%%%%%%%%%%%%%%%%%%%%%%%%%%%%%%%%%%%%%%%%%%%%%%%%%%%%%%%%%%%%%%%%

\section{Introduction}\label{sec:intro}

In this note we generalize the Semple tower from enumerative algebraic geometry (\cite{semple}, \cite{lejeune}, \cite{kennedy}) to study the local theory of Goursat multi-flags, a class of nonholonomic distributions that generalizes the contact distributions of $J(n;k)$, where $J(n;k)$ is the space of $k$-jets of maps from $\mathbb{R}$ to $\mathbb{R}^n$.

A Goursat multi-flag is a nonholonomic distribution $D$ with {\it slow growth}.  More concretely,  a Goursat $n$-flag of length $k$ is a distribution $D$ of rank $(n+1)$ sitting in a 
$(n+1) + kn$ dimensional ambient manifold where the rank of the associated flag of distributions
$$D \hspace{.2in} \subset  \hspace{.2in} D + [D,D]\hspace{.2in} \subset  \hspace{.2in} D + [D,D] + [[D,D],[D,D]]  \hspace{.2in} \subset \dots $$
grows by $n$ at each bracketing step. One often considers a local family of vector fields that locally spans the distribution $D$.
For the case of $n=1$ and $n=2$ there is a concrete Pfaffian description of Goursat multi-flags as the vanishing of locally defined $1$-forms.
Cartan proved that a generic germ of a Goursat $1$-flag can always be described by the vanishing of the $1$-forms
\[
\omega_{1} = dy-u_{1}dx, \, \, \omega_{2} = du_{1} - u_{2}dx, \, \, \cdots  , \, \, \omega_{n} = du_{n-1} - u_{n}dx  
\label{oneform}
\]

\noindent where $(x,y,u_{1}, \cdots , u_{n})$ are local coordinates on $R^{2} \times (S^{1})^{n}$ and showed that any two generic Goursat germs are equivalent to the germ of the distribution described by these $1$-forms (\cite{cartan}).
However, A. Kumpera and C. Ruiz pointed out in \cite{kumpera3} that there are Goursat $1$-flags that are not locally equivalent to the vanishing of the $1$-forms listed above and are referred to as \textit{singularities}.  One way to determine these $1$-forms is by using a coding system known as $RVT$ code developed in \cite{mont2} and \cite{castro}.

In \cite{castro2}, we worked with a tower of fiber bundles, denoted by $S(n;k)$, containing the space $J(n;k)$ as an open dense subset.  Each space $S(n;k)$ comes equipped with a geometric distribution $\Delta_k$.
The link between Goursat multi-flags and the Semple tower is the following:
\begin{quote}
{\it the problem of classifying Goursat multi-flags up to local equivalence is equivalent to the classification of points
in the Semple tower $S(n;k)$ up to symmetry.}
\end{quote}

Our classification is only for Goursat multi-flags of short length. In this note we specialize to the study of $\mbox{Diff}(\mathbb{A}^3)$-orbits in $S(2;4)$. This is equivalent to studying local normal forms of Goursat $2$-flags of length up to $4$.  Our proofs are constructive and rely upon elementary projective geometry and combinatorial invariants of spatial curves (\cite{arnsing}).  For the determination of local orbits we introduce the so called {\it isotropy method}.

%-----------------------------------------------------------------
\section{Definitions and Statement of Main Results}

A {\it geometric distribution} hereafter denotes a linear subbundle of the tangent bundle.

\subsection{Prolongation}
Let the pair $(Z,\Delta)$ denote a manifold $Z$ of dimension $d$ equipped with a distribution $\Delta$ of rank $r$.  We denote by $\Pj (\Delta)$ the {\it projectivization} of $\Delta$.
As a manifold, $$\Pj (\Delta) := Z^1,$$ has dimension $d + (r-1)$.

Given an analytic horizontal curve $c:(I,0)\rightarrow (Z,q)$, where $I$ is some open interval in $\mathbb{R}$ containing the origin and $c(0)=q$, we can naturally define a new curve $$c^{1}:(I,0)\rightarrow (Z^1,(q,\ell))$$ with image in $Z^1$ and where $\ell = \mbox{span}\{ \frac{dc}{dt}(0) \} \subset \Delta_{q}$.
This new curve, $c^{1}(t)$, is called the \textit{prolongation} of $c(t)$.  This procedure can be iterated.

The manifold $Z^1$ also comes equipped with a distribution $\Delta_1$ called the {\it Cartan prolongation of $\Delta$} (\cite{bryant}) which is defined as follows.
Let $\pi : Z^1 \rightarrow Z$ be the projection map $(p, \ell)\mapsto p$. Then
$$\Delta_1(p,\ell) = d\pi_{(p,\ell)}^{-1}(\ell),$$
It is easy to check that $\Delta_1$ is also a  distribution of rank $r$.

\begin{table}[t]
\bitwonumcaption{Short}{Some geometric objects and their Cartan prolongations.}{Tableau}{Court}{Certains objets g\'eom\'etriques et leurs prolongements de Cartan.}
\begin{tabular}{|c|c|}
  \hline
   curve $c:(I,0)\rightarrow (Z,q)$ & curve $c^{1}:(I,0)\rightarrow (Z^1,q),$ \\
   & $c^{1}(t) =\text{(point,moving line)} = (c(t),\mbox{span}\{ \frac{dc}{dt}(t) \}) $\\ \hline
  diffeomorphism $\Phi: Z \circlearrowleft$ & diffeomorphism $\Phi^{1}: Z^1 \circlearrowleft$, \\
  & $\Phi^{1}(p,\ell) = (\Phi(p),d\Phi_p(\ell))$ \\ \hline
  rank $r$ linear subbundle  & rank $r$ linear subbundle $\Delta_{1 (p, \ell)}=d\pi_{(p,\ell)}^{-1}(\ell)\subset TZ^1$,\\
  $\Delta \subset TZ$ & $\pi: Z^1 \rightarrow Z$ is the canonical projection. \\
  \hline
\end{tabular}
\label{tab:prolong}
\end{table}

By a {\it symmetry} of the pair $(Z,\Delta)$ we mean a local diffeomorphism $\Phi$ of $Z$ that preserves the subbundle $\Delta$.

The symmetries of $(Z,\Delta)$ can also be prolonged to symmetries $\Phi^{1}$ of $(Z^1,\Delta_{1})$ as follows.
Define
$$\Phi^{1}(p,\ell):=(\Phi(p), d\Phi_{p}(\ell)).$$
Since\footnote{We also use the notation $\Phi_{*}$ for the pushforward or tangent map $d\Phi$.}  $d\Phi_{p}$ is invertible, the second component is well defined as a projective map. This new transformation in $(Z^1,\Delta_1)$ is the \textit{prolongation} of $\Phi$.
Objects of interest and their Cartan prolongations are summarized in Table $1$. 
%\ref{tab:prolong}

\subsection{Constructing the Semple Tower.}
\begin{defn}
The \textit{Semple Tower} is a sequence of manifolds with distributions, $(S(n;k), \Delta_{k})$,
together with fibrations $$\cdots \rightarrow S(n;k) \rightarrow S(n;k-1) \rightarrow \cdots \rightarrow
S(n;1) \rightarrow S(n;0) = \mathbb{A}^{n+1}$$
and we write $\pi_{k,i}: S(n;k) \rightarrow S(n;k)$ for the respective bundle projections.
\end{defn}

\begin{defn}
$\mbox{ Diff}(n)$ is taken to be the \textit{pseudogroup of diffeomorphism germs of} $\mathbb{A}^{n}.$
\end{defn}

\begin{thm}
\label{thm:sym}
For $n > 1$ and $k>0$ any local diffeomorphism of $S(n;k)$ preserving the distribution $\Delta_{k}$
is the restriction of the $k$-th prolongation of a local diffeomorphism $\Phi \in \mbox{Diff}(n)$.
\end{thm}

Proof: See \cite{yamaguchi1}, pg. $795$.  $\square$

Shibuya and Yamaguchi also point out that this is a result due to A. B\"{a}cklund (\cite{backlund}).

\begin{rem}
From now on we will write $S(k)$ to denote $S(2;k)$.
\end{rem}

\begin{defn}
Two points $p,q$ in $S(k)$ are said to be \textit{equivalent}, written $p \sim q$, if there is a $\Phi\in \mbox{Diff}(3)$ such that $\Phi^{k}(p)=q$.
\end{defn}

\begin{defn}
Let $p \in S(k)$ then we denote by $\mathcal{O}(p)$ to be the  {\it orbit of the point $p$} under
the action by elements of $\mbox{\em Diff}(3)$ on the $k$-th level of the Semple Tower.  This means that
a point $q$ is an element in $\mathcal{O}(p)$ if $q$ is equivalent to the point $p$.
\end{defn}

 \begin{figure}
 \raggedright
    \subfloat[Prolongation of a distribution]{\label{fig:prolongdist} \def\svgwidth{.3 \columnwidth}
     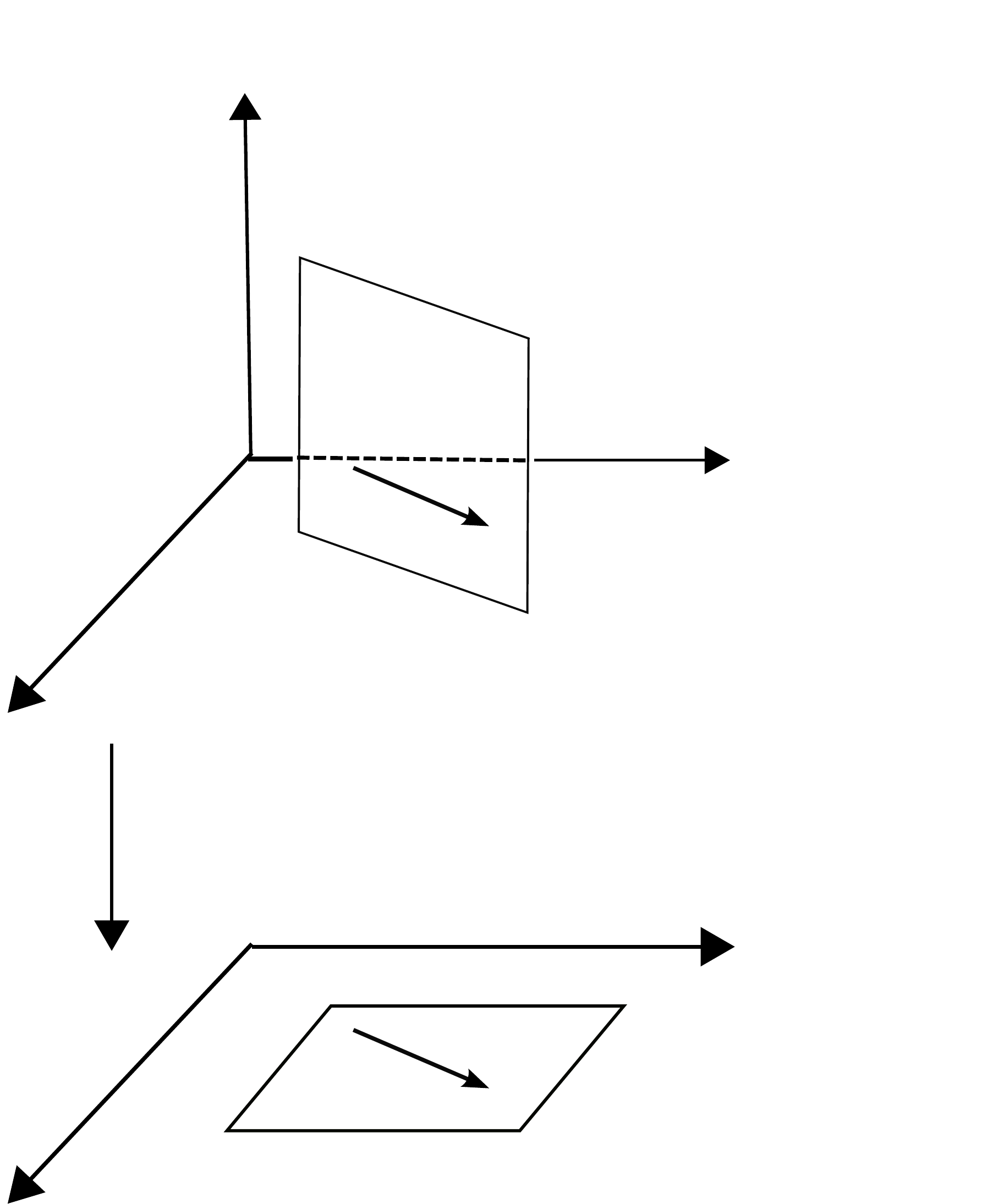}
    \subfloat[Prolongation of a diffeomorphism]{\label{fig:prolongdiff} \def\svgwidth{.3 \columnwidth}
     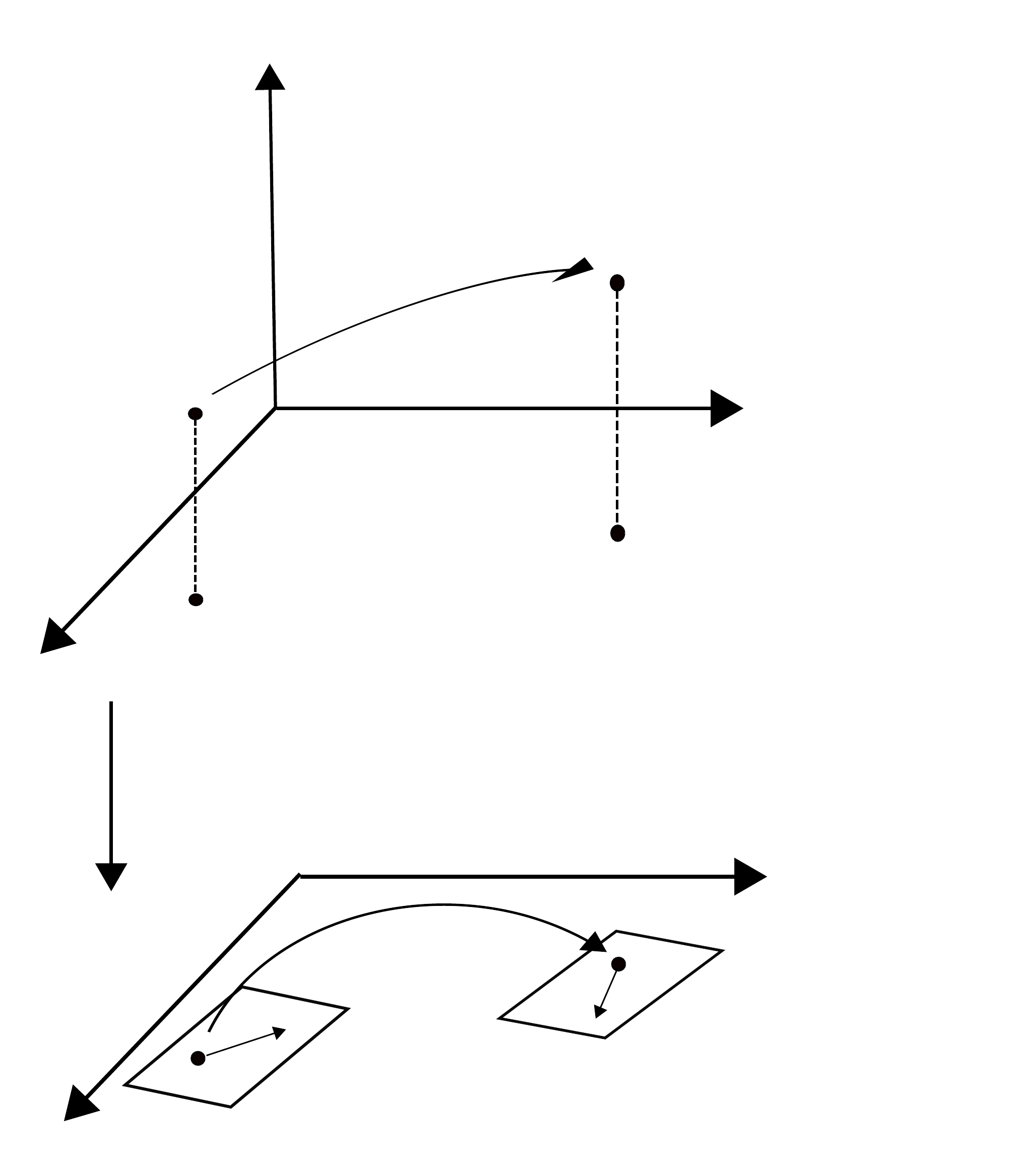}
     \subfloat[Prolongation of a curve]{\label{fig:prolongcur} \def\svgwidth{.3 \columnwidth}
     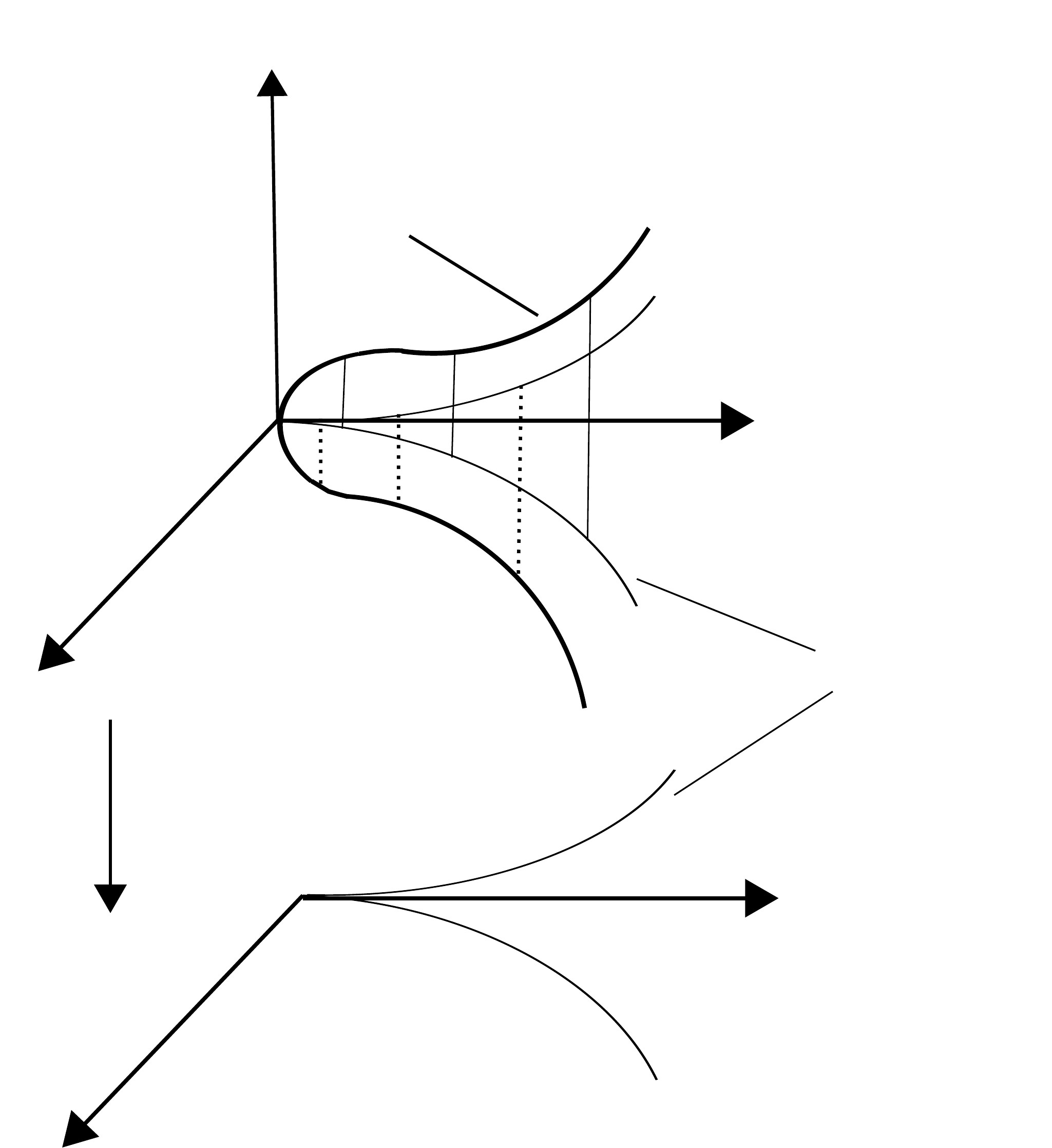}
     %%\caption{Prolongations}
     \bitwonumcaption{Short}{Prolongations.}{Figure}{Court}{Prolongements.}
     \label{fig:prolongation}
  \end{figure}

\begin{defn}
We say that a curve or curve germ $\gamma: (\R, 0) \rightarrow (\mathbb{A}^{3}, p_{0})$ \textit{realizes} the point $p_{k} \in S(k)$
if $\gamma^{k}(0) = p_{k}$, where $p_{0} = \pi_{k,0}(p_{k})$.
\end{defn}

% ----------------------------------------------------------------
\section{Main Results and Ideas of Proofs}

From our recent work in \cite{castro2} we have been able to completely classify the points within the first four
levels of the tower.  Our findings are summarized by the following.

\begin{thm}[Orbit counting per level]\label{thm:main}
In the $n=2$ (or spatial) Semple Tower the number of orbits within each of the first four levels of the tower are as follows:
\begin{itemize}
\item Level $1$ has $1$ orbit,
\item Level $2$ has $2$ orbits,
\item Level $3$ has $7$ orbits,
\item Level $4$ has $34$ orbits.
\end{itemize}
\end{thm}

\

In order to classify each of these orbits we used a coding system known as $RVT$ coding which partitions the points
into equivalence classes that are invariant under the action of the symmetries at any level.  The details
of this coding system can be found in \cite{castro} and \cite{castro2}.

\begin{thm}[Listing of orbits within each $RVT$ code]\label{thm:count}
Table $2$ is a breakdown of the number of orbits that appear within each $RVT$ class within the first three levels.

\begin{table}[t!]
\bitwonumcaption{Short}{Number of orbits within the first three levels of the Semple Tower.}{Tableau}{Court}{Nombre d'orbites dans les trois premiers niveaux de la tour de Semple.}
\begin{tabular}{|c|c|c|c|}
\hline
Level of tower &      $RVT$ code  & Number of orbits & Normal forms \\
$1$            &        $R$       & $1$              &  $(t,0,0)$ \\
$2$            &        $RR$      & $1$              &  $(t,0,0)$ \\
               &        $RV$      &  $1$             &  $(t^{2}, t^{3}, 0)$\\
$3$            &       $RRR$      & $1$              &  $(t,0,0)$ \\
               &       $RRV$      & $1$              &  $(t^{2}, t^{5}, 0)$ \\
               &       $RVR$      & $1$              &  $(t^{2}, t^{3}, 0)$ \\
               &       $RVV$      & $1$              &  $(t^{3}, t^{5}, t^{7})$, $(t^{3}, t^{5}, 0)$ \\
               &       $RVT$      & $2$              &  $(t^{3}, t^{4},t^{5})$, $(t^{3}, t^{4}, 0)$ \\
               &       $RVL$      & $1$              &  $(t^{4}, t^{6}, t^{7})$ \\
\hline
\end{tabular}
\label{tab:codes}
\end{table}

\noindent For level $4$ there is a total of $23$ possible $RVT$ classes.  Of the $23$ possibilities $14$ of them consist of a single orbit.
The classes $RRVT$, $RVRV$, $RVVR$, $RVVV$, $RVVT$, $RVTR$, $RVTV$, and $RVTL$ consist of $2$ orbits, and the class $RVTT$ consists of $4$ orbits.
\end{thm}

\begin{rem}
The curve normal forms were computed by turning $\mathbb{A} ^{3}$ into the vector space $\R ^{3}$ and doing our computations over the origin in $\R ^{3}$. 
More details about the computations and properties of the curve normal forms can be found in \cite{castro2}.
\end{rem}

One recent development that was not included in \cite{castro2} is the following.

\begin{thm}
\label{thm:new}
Let $\omega$ be an $RVT$ class that has a total of $k$ orbits.
\

\noindent Then the addition of $R$'s to the beginning of the code $\omega$, meaning $R \cdots R \omega$, will be an $RVT$ class with $k$-orbits.
\end{thm}

\noindent Theorem \ref{thm:new} allows one to reduce the number of calculations needed to compute the number of orbits within an
$RVT$ class found within higher levels of the Semple Tower.

\subsection{Idea of Proofs}

When we began our classification of the points within the Semple Tower we followed the approach taken by Montgomery and Zhitomirskii in \cite{mont2} by using the singular curves approach as well as the $RVT$ coding system, which partitioned the various points within the tower (\cite{mont2}, \cite{castro}).  The curve approach yielded the complete classification of the points of the plane Semple Tower $S(1;k)$.  However, while classifying the points within level $4$ of the tower we noticed that our curve approach broke down very quickly.  This lead us to develop the \textit{isotropy method}.  While this technique does not yield the same normal forms as the curve approach, it did help us gain a better overall picture of some of the subtleties about the Semples Tower. In particular, this method helped us understand the incidence relations between the various critical hyperplanes and how these critical directions arise within the tower (\cite{castro}, \cite{castro2}).  One other important feature of the isotropy method is that it can be used to classify the points within any $RVT$ class within the Semple Tower.

\subsection{The isotropy method.}
\

Suppose we want to look at a particular $RVT$ class, at the $k$-th level, given by $\omega$ (a word of
length $k$) and we want to see how many orbits there are.  Suppose as well that we understand its projection
$\pi_{k,k-1}(\omega)$ one level down, which decomposes into $N$ orbits.  Choose representative points $p_{i}$, $i = 1, \cdots , N$ for
the $N$ orbits in $\pi_{k,k-1}(\omega)$, and consider the group $G_{k-1}(p_{i})$ of level $k-1$ symmetries that fix $p_{i}$.  This group is called the \textit{isotropy group of} $p_{i}$.  Since elements $\Phi^{k-1}$ of the isotropy group
fix $p_{i}$, their prolongations $\Phi^{k} = (\Phi^{k-1}, \Phi^{k-1}_{\ast})$ act on the fiber over $p_{i}$.
Under the action of the isotropy group the fiber decomposes into some
number $n_{i} \geq 1$ (possibly infinite) of orbits.  Summing up, we find that $\omega$ decomposes into
$\sum_{i = 1}^{N} n_{i} \geq N$ orbits.
This will tell us how many orbits there are for the class $\omega$.

\begin{figure}
 \def\svgwidth{.4 \columnwidth}
 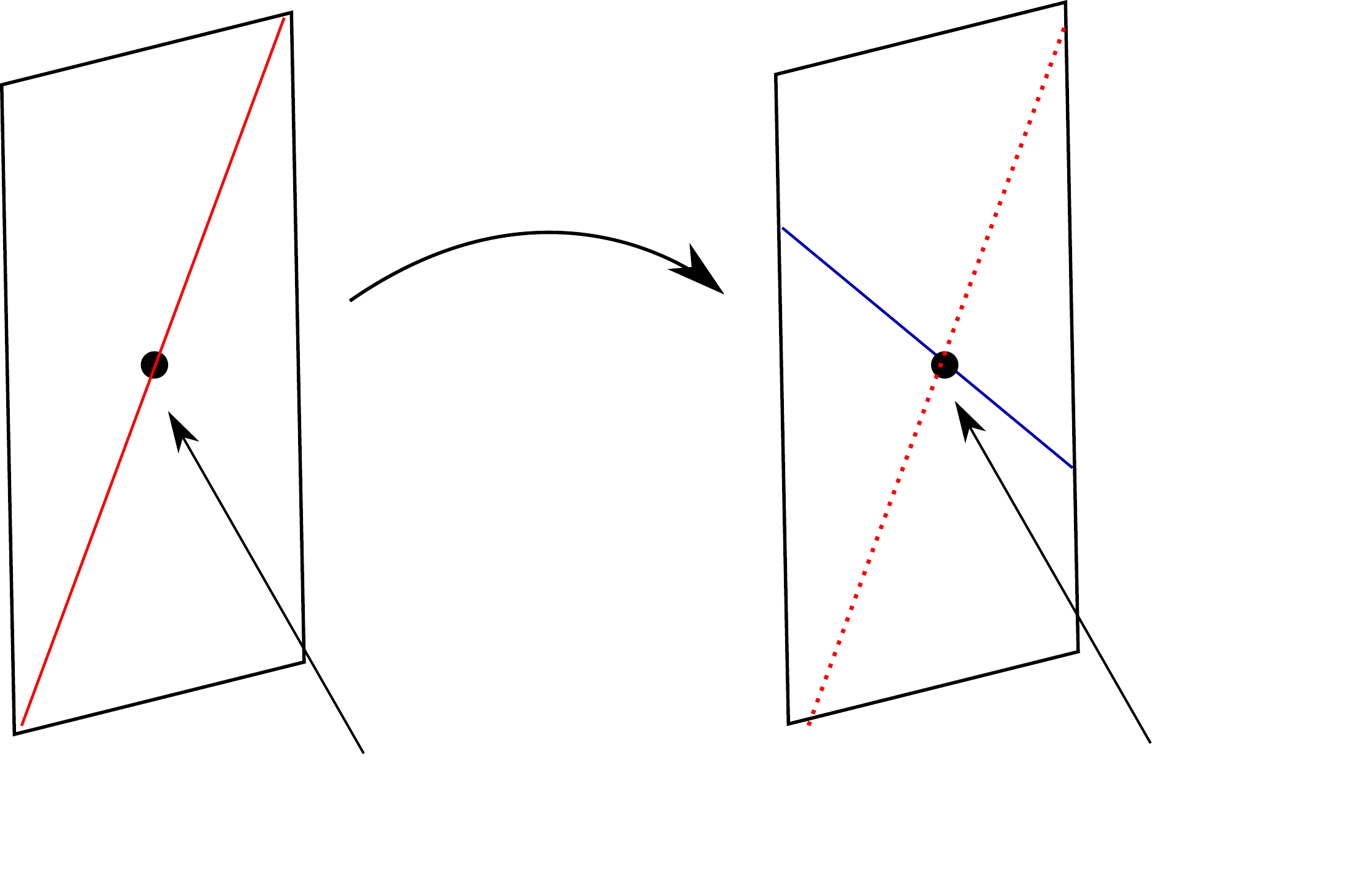
\bitwonumcaption{Short}{The isotropy method.}{Figure}{Court}{Utilisation de l'isotropie de l'action de groupe pour compter orbites.}
 \label{fig:isotropy}
\end{figure}

This is the idea behind the approach (see Figure $2$).  We want to point out that Theorem \ref{thm:sym} is a crucial result which allows our isotropy method to work.

\begin{rem}[The class $RVVV$]
In our work we carried out the explicit details of how one uses the isotropy method to calculate the number of orbits.  In particular, we compute the orbits of the $RVT$ class $RVVV$ in explicit detail using the isotropy method.  One finds there is a total of $2$ orbits within the class
$RVVV$.
\end{rem}

% ----------------------------------------------------------------
\section{Conclusion}
\

Our work with the Semple Tower has been primarily focussed on classifying points up to symmetry at various levels within the tower.  As was pointed out in \cite{mont2}, it is possible to end up with moduli appearing within various $RVT$ classes, meaning that there is a whole continuum's worth of points in a particular class that are inequivalent.  We are interested in determining at which level does this first occur in the $n=2$ Semple Tower.  We believe that one starting point is to look at the $RVT$ classes of the form $RV^{k}R^{l}VR^{m}$ for $k,l,m$ positive integers.
\

In a series of papers by J.C. Hausmann and E. Rodriguez (\cite{rodriguez1}, \cite{hausmann}) they examine a control theory problem known as the \textit{Snake Charmer Algorithm}.  This is the continuous analogue of the polygonal snakes that were initially studied by Hausmann in \cite{hausmann}.  Due to the fact that the configuration space for snakes involves looking at continuous $C^{1}$-curves on $S^{n}$, Hausmann and Rodriguez take the group action on the configuration space to be the M\"obius group.  This problem also shares a connection to the world of Goursat-multi flags.  Pelletier and Slayman in \cite{pelletier2}  show that Hausmann's polygonal snakes correspond to the \textit{$k$-bar system}  and establishes a connection between the distribution arising from the configuration space of snakes to the theory of Goursat multi-flags (\cite{pelletier}).  We are interested in extending their work using geometric mechanics to understand what sort of geometric phenomena can appear when trying to control a polygonal snake.

%\section{}
%\subsection{}

% -----------------------------------------------------------------

%%%%%%%%%%%%%%%%%%%%%%%%%%%%%%%%%%%%%%%%%%%%
\medskip
% \bibliographystyle{alpha}  %alpha style uses authors last name and year
% \bibliography{orbits-references}

%%%%%%%%%%%%%%%%%%%%%%%%%%%%%%%%%%%%%%%%%%%%%

\end{document}